%% file: psiuni.tex
\documentclass[11pt]{article}

\usepackage{amsmath}
\usepackage{amsfonts}
\usepackage{amssymb}
\usepackage{amsthm}
\usepackage{graphics}
\usepackage{amscd}
\usepackage{graphicx,epsfig}
\usepackage{color}
\usepackage[all]{xy}
\usepackage{url}

\DeclareMathOperator{\Div}{div}

\renewcommand{\epsilon}{\varepsilon}

\newcommand{\R}{\mathbb{R}}

\renewcommand{\L}{\mathbb{L}}

\newcommand{\N}{\mathbb{N}}
\newcommand{\la}{\langle}
\newcommand{\ra}{\rangle}
\newcommand{\eps}{\varepsilon}

\newcommand{\dd}{\mathrm{d}}

\newcommand{\Ome}{\Omega}
\newcommand{\Omet}{\widetilde{\Omega}}

\newtheorem{thm}{Theorem}
\newtheorem{cor}[thm]{Corollary}

\newtheorem{lem}[thm]{Lemma}

\renewcommand{\phi}{\varphi}
\newcommand{\dis}{\displaystyle}

\newtheorem*{thm*}{Theorem}

\theoremstyle{remark}

\newtheorem*{rem*}{Remark}

\newcounter{remark}

\newcounter{case}

\newcounter{construction}

\newcounter{fact}

\title{A uniqueness result for maximal surfaces in Minkowski $3$-space}
\author{Laurent Mazet}
\date{}

\begin{document}
\maketitle

\begin{abstract}
In this paper, we study the Dirichlet problem associated to the
maximal surface equation.  We prove the uniqueness of bounded
solutions to this problem in unbounded domain in $\R^2$.
\end{abstract}

\section*{Introduction}
\input psiuni0.tex

\section{The uniqueness result}
\input psiuni1.tex

\end{document}

%% file: psiuni0.tex
We consider the Minkowski space-time $\L^3$ \textit{i.e.} $\R^3$ with
the following pseudo-euclidean metric $\la x,y \ra= x_1y_1+ x_2y_2-
x_3y_3$. We define $|x|_{\L}^2=\la x,x \ra$.

A vector is said to be \emph{spacelike} if $|x|_\L^2>0$ and a surface
$S$ of class $C^1$ is said to be spacelike if $|\cdot |_\L^2$ is
positive definite on the tangent space to $S$. Such a surface is
locally the graph of a function over a domain in $\R^2$.

If $v$ is a function in a domain $\Ome$ in $\R^2$ (in the paper, we
always assume that $\Ome$ has smooth boundary), the graph of $v$ is
spacelike if and only if $|\nabla v|<1$. The function $v$ is then Lipschitz
continuous and it extends to the closure $\overline{\Ome}$. We denote
by $\phi$ the trace $v|_{\partial\Ome}$ of $v$ on the boundary. 

The maximal area problem in the class of spacelike surfaces consists
in solving the following variational problem:
$$\max_v \int_\Ome \sqrt{1-|\nabla v|^2}\dd x,\quad v|_{\partial\Ome}=
\phi$$ 

The critical points of this functional are the solutions of the maximal
surface equation :
\begin{equation*}
\Div\frac{\nabla v}{\sqrt{1-|\nabla v|^2}}=0
\tag{$*$}
\end{equation*}
The maximal area problem is then linked to the Dirichlet problem
associated to \eqref{Mse}: to find a solution $v$ of \eqref{Mse} in
$\Ome$ such that $v|_{\partial\Ome}=\phi$. This Dirichlet problem has
been already studied by several authors, for exemple see \cite{BS} and
\cite{KM}. 

In this paper, we prove the uniqueness of bounded solutions to the
Dirichlet problem. More precisely, if $\Ome$ is an unbounded domain
and $\phi$ is a bounded continuous function on $\partial \Ome$, we
prove that, if it exists, a solution $v$ of \eqref{Mse} in $\Ome$ with
$v|_{\partial\Ome}=\phi$ is unique (Theorem \ref{main}). The study of
the uniqueness is important in the construction of certain moduli
spaces of maximal surfaces (see \cite{FLS1}and \cite{FLS2}).

This uniqueness result for the maximal surface equation is also
important for the study of the Dirichlet problem associated to the
minimal surface equation. The graph of a function $u$ over a domain 
$\Ome\subset\R^2$ is a surface in $\R^3$ with its standard euclidean
metric, it has vanishing mean curvature if $u$ satisfies the following 
partial differential equation:
\begin{equation*}\label{mse}
\Div\frac{\nabla u}{\sqrt{1+|\nabla u|^2}}=0
\tag{$**$}
\end{equation*} 
This equation implies that there exists  locally  a function $v$ such
that:
$$
\dd v=\dd \Psi_u=\frac{u_x}{\sqrt{1+|\nabla u|^2}}\dd y-
\frac{u_y}{\sqrt{1+|\nabla u|^2}} \dd x
$$
(here $u_x$ and $u_y$ are the first derivatives of $u$). $v=\Psi_u$ is
called the conjugate function to $u$ and a simple computation
shows that $v$ is a solution of \eqref{Mse}. Then the uniqueness for
solutions of \eqref{Mse} can implies uniqueness for solutions of
\eqref{mse}. 

The proof of our uniqueness result uses the same technic as P.~Collin and
R.~Krust in \cite{CK}. But to apply this technic, we need to prove an
estimate for the first derivatives of $v$  in a subdomain of $\Ome$;
this is what we do in Lemma \ref{tech2}.

%% file: psiuni1.tex
Let $\Ome\subset \R^2$ be a domain and $v$ a solution of the maximal
surface equation :
\begin{equation*}\label{Mse}
\Div\frac{\nabla v}{\sqrt{1-|\nabla v|^2}}=0
\tag{$*$}
\end{equation*}
In the following, the quantity $\sqrt{1-|\nabla v|^2}$ will be denoted
by $w_v$. We then define the $1$-form $\alpha_v$ as follows: 
$$\alpha_v=\frac{v_x}{w_v}\dd y-\frac{v_y}{w_v}\dd x$$
where $v_x$ and $v_y$ are the first derivatives of $v$. The maximal
surface equation is then equivalent to $\dd \alpha_v=0$.

First, we need a technical lemma.
\begin{lem}\label{tech}
Let $v$ and $v'$ be two functions. Let $P$ be a point
in $\Ome$ and $\eps>0$ such that $|\nabla v|(P)\le 1-\eps$ and
$|\nabla v'|(P)\le 1-\eps$. Then there exists a constant $C(\eps)$ that
depends only on $\eps$ such that, at the point $P$, we have:
\begin{equation}\label{lem}
\left( (\nabla v-\nabla v')\cdot \left(\frac{\nabla
v}{w_v}-\frac{\nabla v'}{w_{v'}}\right)\right) \ge
C(\eps)\left|\frac{\nabla v}{w_v}-\frac{\nabla v'}{w_{v'}}\right|^2  
\end{equation}
\end{lem}
\begin{proof}
We define $n=(-v_x,-v_y,1)/w_v$ and $n'=(-v'_x,-v'_y,1)/w_{v'}$. We
have $|n|_\L^2=-1$ and $|n'|_\L^2=-1$, then 
\begin{align*}
\left( (\nabla v-\nabla v')\cdot \left(\frac{\nabla
v}{w_v}-\frac{\nabla v'}{w_{v'}}\right)\right)&=  \la
(w_{v'}n'-w_vn),(n'-n)\ra\\ 
&=(w_v+w_{v'})(-1-\la n,n'\ra)\\
&=\frac{w_v+w_{v'}}{2} |(n'-n)|_\L^2
\end{align*}
Since $|\nabla v|\le 1-\eps$ and $|\nabla v'|\le 1-\eps$ there
exists $C_1(\eps)>0$ such that 
\begin{equation}\label{eq1}
(w_v+w_{v'})/2\ge C_1(\eps)
\end{equation}
Besides
$$\dis|(n'-n)|_\L^2 =\left|\frac{\nabla v}{w_v}-\frac{\nabla
v'}{w_{v'}}\right|^2  - \left(\frac{1}{w_v}- \frac{1}{w_{v'}}\right)^2$$

Let $x\in\R^2$ be $\nabla v/w_v$ and $x'$ be
  $\nabla v'/w_{v'}$. Then $1/w_v=\sqrt{1+|x|^2}$ and
  $1/w_{v'}=\sqrt{1+|x'|^2}$. Since $\nabla v$ and $\nabla v'$ are
  bounded by $1-\eps$, there exists $R(\eps)$ such that $|x|$ and
  $|x'|$ are bounded by $R(\eps)$. Then: 
\begin{align*}
\frac{|(n'-n)|_\L^2}{\left|\frac{\nabla v}{w_v}-\frac{\nabla
  v'}{w_{v'}}\right|^2}&= 1- \frac{\left(\frac{1}{w_v}-
  \frac{1}{w_{v'}}\right)^2}{|x-x'|^2}\\
&=1-\frac{\left(\sqrt{1+|x|^2}-\sqrt{1+|x'|^2}\right)^2}{|x-x'|^2}\\
&=1-\frac{\left(|x|^2-|x'|^2\right)^2} {|x-x'|^2\left(\sqrt{1+|x|^2}+
  \sqrt{1+|x'|^2}\right)^2}\\
&=1-\left(\frac{|x|-|x'|}{|x-x'|}\right)^2
  \left(\frac{|x|+|x'|}{\left( \sqrt{1+|x|^2}+
  \sqrt{1+|x'|^2}\right)}\right)^2 \\
&\ge 1-\left(\frac{|x|+|x'|}{\left( \sqrt{1+|x|^2}+
  \sqrt{1+|x'|^2}\right)}\right)^2 >0  
\end{align*}
By continuity and since $|x|$ and $|x'|$ are bounded by $R(\eps)$,
there exists a constant $C_2(\eps)>0$ such that:
\begin{equation}\label{eq2}
1-\left(\frac{|x|+|x'|}{\left( \sqrt{1+|x|^2}+
\sqrt{1+|x'|^2}\right)}\right)^2 >C_2(\eps)
\end{equation}
Then in combining \eqref{eq1} and \eqref{eq2}, we get \eqref{lem} with
$C(\eps)=C_1(\eps)C_2(\eps)$. 
\end{proof}

We denote by $d$ the usual distance in $\R^2$ and by $d_\Ome$ the
intrinsic metric in $\Ome$ \textit{i.e.} $d_\Ome(p,q)$ is the infimum
of the length of all paths in $\Ome$ going from $p$ to $q$. Let
$\delta>0$, we denote by $\Ome_\delta$ the set $\{p\in\Ome\,|\,
d_\Ome(p,\partial\Ome)>\delta\}$. We then can write our uniqueness
result.

\begin{thm}\label{main}
Let $\Ome$ be an unbounded domain in $\R^2$ and $\phi$ a bounded
continuous function on $\partial\Ome$. Let $v$ and $v'$ be two bounded
solutions of \eqref{Mse} in $\Ome$ with
$v|_{\partial\Ome}=\phi=v'|_{\partial\Ome}$. Then $v=v'$. 
\end{thm}

\begin{proof}
Let $v$ and $v'$ be two such solutions. We assume that $\sup v-v'>0$
and we denote this supremum by 
$4\delta$. Let $a\in [2\delta,3\delta]$ be chosen such that
$\Omet=\{v>v'+a\}$ has smooth boundary. Since $2\delta \le a \le
3\delta$ and $v$ and $v'$ are $1$-Lipschitz continuous $\Omet\subset
\Ome_\delta$. We then have the following lemma.
\begin{lem}\label{tech2}
There exists $\eps>0$ such that, in
$\Omet$, $|\nabla v|\le 1-\eps$ and $|\nabla v'|\le 1-\eps$.
\end{lem}

Before proving this lemma, we finish Theorem \ref{main} proof. Let
$\tilde v$ denote $v-v'-a$ and $\tilde\alpha$ denote
$\alpha_v-\alpha_{v'}$. 

For $r>0$, we define $\Omet_r=\{p\in\Omet\,|\, |p|<r\}$ and
$C_r=\{p\in\Omet\,|\, |p|=r\}$. Since $\tilde v=0$ on $\partial
\Omet_r\backslash C_r$ and $\tilde\alpha$ is closed, we have :
$$\int_{C_r}\tilde v\tilde\alpha=\int_{\partial\Omet_r}\tilde
v\tilde\alpha= \iint_{\Omet_r}\dd\tilde v\wedge\tilde\alpha$$
Since $\dd\tilde v\wedge\tilde\alpha =  \left( (\nabla v-\nabla
v')\cdot \left(\frac{\nabla v}{w_v}-\frac{\nabla
v'}{w_{v'}}\right)\right) \dd x\wedge \dd y$, Lemma \ref{tech} and
Lemma \ref{tech2} imply that:
$$C(\eps)\iint_{\Omet_r}|\tilde\alpha|^2\le \int_{C_r}\tilde
v\tilde\alpha$$

Let $r_0$ be such that $\dis\mu= C(\eps)\iint_{\Omet_{r_0}}
|\tilde\alpha|^2>0$. In $\Omet$, $\tilde v$ is bounded by $2\delta$
then :
$$\mu+C(\eps)\iint_{\Omet_r\backslash\Omet_{r_0}}|\tilde\alpha|^2\le
2\delta\int_{C_r}|\tilde\alpha|$$

Let us denote $\int_{C_r}|\tilde\alpha|$ by $\eta(r)$. Then by
Schwartz's Lemma :
$$
\eta^2(r)\le \ell(C_r)\int_{C_r} |\tilde\alpha|^2\le 2\pi r \int_{C_r}
|\tilde\alpha|^2 
$$
Then $\dis \int_{C_r} |\tilde\alpha|^2 \ge \frac{\eta^2(r)}{2\pi r}$
and
$$\int_{r_0}^r \frac{\eta^2(t)}{2\pi t}\le \iint_{\Omet_r\backslash
\Omet_{r_0}} |\tilde\alpha|^2$$
Then :
\begin{equation}\label{equadiff}
\mu+C(\eps)\int_{r_0}^r \frac{\eta^2(t)}{2\pi t}\le 2\delta\eta(r)
\end{equation}

Let $y$ be the solution of the following Cauchy problem :
$$ y'(t)=C(\eps)\frac{y^2(t)}{4\pi\delta t},\quad y(r_0)=
\frac{\mu}{4\delta}$$
$y$ is defined on $[r_0,r_1)$ with
$r_1=r_0\exp(\frac{16\pi\delta^2}{\mu C(\eps)})$ and is defined by :
$$\frac{4\delta}{\mu}-\frac{1}{y(t)}=\frac{C(\eps)}{4\pi \delta}
\ln\frac{t}{r_0}$$ 

By \eqref{equadiff}, $\eta(t)\ge y(t)$ on $[r_0,r_1)$ and, since
$\dis\lim_{t\rightarrow r_1}y(t)=+\infty$, we get a
contradiction, indeed $\eta$ is continuous . Then $v=v'$. 
\end{proof}

As we say in the introduction Theorem \ref{main} has a consequence for
solution of the minimal surface equation.

\begin{cor}
Let $\Ome$ be an unbounded simply-connected domain in $\R^2$. Let $u$
and $u'$ be two solutions of \eqref{mse} in $\Ome$ such that  $\Psi_u$
and $\Psi_{u'}$ are bounded in $\Ome$ and $\Psi_u=\Psi_{u'}$ on
$\partial \Ome$. Then $u-u'$ is constant.
\end{cor}

We need the simple-connecteness hypothesis to ensure that $\Psi_u$ and
$\Psi_{u'}$ are well defined. 

\begin{proof}
$\Psi_u$ and $\Psi_{u'}$ are two solutions of \eqref{Mse} in $\Ome$,
then, by theorem \ref{main}, $\Psi_u=\Psi_{u'}$. Then $\nabla u=\nabla
u'$ and $u-u'$ is constant.
\end{proof}

To end Theorem \ref{main} proof, we have to prove Lemma~\ref{tech2}.  

\section{The gradient estimate}

This section is devoted to the proof of the gradient estimate in Lemma
\ref{tech2}; This is the last step in Theorem \ref{main} proof.

\begin{proof}[Proof of Lemma \ref{tech2}]
If Lemma \ref{tech2} is not true, we can assume that
$\sup_{\Omet}|\nabla v|=1$. Then there exists 
$(p_n)$ a sequence in $\Omet$ such that $|\nabla v|(p_n)\rightarrow
1$. Let $O$ be the point $(0,0)$. Let $r_n$ be the affine rotation in
$\R^2$ such that $r_n(O)=p_n$ and $R_n^{-1}\big(\nabla
v(p_n)\big)=(|\nabla v|(p_n),0)$ 
($R_n$ is the linear rotation 
associated to $r_n$). We then define $v_n=v\circ r_n$ which is a
solution of \eqref{Mse} in $\Ome_n=r_n^{-1}\Ome$. We have $\nabla
v_n=R_n^{-1}\nabla v$ then $\nabla v_n(O)\rightarrow (1,0)$. In the
same way we define $v'_n=v'\circ r_n$.

Let $I(a,b)\subset\R^2$ be the segment $[a,b]\times\{0\}$ ($a<b$). Let 
$\eps$ be positive, $\eps$ will be fixed later but let us notice that
$\eps/\delta$ will be small. Let $D(a,b)$ denote the set $\{p\in\R^2\,
|\, d(p,I(a,b))<\eps\}$, $D(a,b)$ is the union of a rectangle of width
$2\eps$ and length $b-a$ and two half-disks of radius $\eps$. 

For every $n$, we define $a_n$ and $b_n$ by: $a_n=\inf\{a\le0\,|\,
D(a,0)\subset \Ome_n\}$ and $b_n=\sup\{b\ge0\,|\, D(0,b)\subset
\Ome_n\}$. Since $\eps<\delta$ and $O\in {\Ome_n}_\delta$ (because
$p_n\in\Ome_\delta$), $b_n>0$ and $a_n<0$; besides $D(a_n,b_n)\subset
\Ome_n$. We define $b_\infty=\liminf b_n$, $b_\infty>0$, $b_\infty$
may take the value $+\infty$; by taking a subsequence, we assume that
$b_\infty=\lim b_n$. Then we define $a_\infty=\limsup a_n$,
$a_\infty<0$, $a_\infty$ may take the value $-\infty$; as above we can
assume that $a_\infty=\lim a_n$. Let $\beta\le\min(\eps/2,
|a_\infty|,b_\infty)$, let A denote $a_\infty +\beta$ if
$a_\infty>-\infty$ and any negative number if not and $B$ denote
$b_\infty-\beta$ if $b_\infty<+\infty$ and any positive number if
not. For $n$ big enough, $D(A,B)\subset \Ome_n$ (see  Figure
\ref{dessin}). 

\begin{figure}
\begin{center}
\resizebox{0.8\linewidth}{!}{\input{dessin.pstex_t}}
\caption{\label{dessin}}
\end{center}
\end{figure}

Since  $D(A,B)$ is simply connected, for each big $n$ in $\N$, there
exists $u_n$ a function on $D(A,B)$ such that $\dd
u_n=\alpha_{v_n}$. Besides the  function $u_n$ satisfies the minimal
surface equation: 
\begin{equation*}
\Div\frac{\nabla u_n}{\sqrt{1+|\nabla u_n|^2}}=0
\tag{\ref{mse}}
\end{equation*}
The graph of $u_n$ is a minimal surface in $\R^3$ with the euclidean
metric. We have
$$\dd v_n=\frac{{u_n}_y}{\sqrt{1+|\nabla u_n|^2}}\dd x-
\frac{{u_n}_x}{\sqrt{1+|\nabla u_n|^2}}\dd y$$
Then $v_n$ is the opposite of the conjugate function to $u_n$. Since
$\nabla v_n(O)\rightarrow (1,0)$, $|\nabla u_n|(O)\rightarrow +\infty$
and $\frac{\nabla u_n}{|\nabla u_n|}(O)\rightarrow (0,1)$. Then
$\{y=0\}\cap D(A,B)$ is a line of divergence
for the sequence $(u_n)$ (see \cite{Ma1,Ma2}). This implies that if
$A-\eps<s <t <B+\eps$:
\begin{equation}\label{truc}
\lim v_n(t,0)-v_n(s,0)=t-s
\end{equation}

By hypothesis, $v$ is bounded by one $M>0$ then $v_n$ is bounded by
$M$. This implies that $A$ and $B$ are bounded thus $a_\infty$ and
$b_\infty$ can not take infinite value; indeed \eqref{truc} implies
$B-A\le 2M$. Then $A=a_\infty+\beta$ and $B=b_\infty-\beta$. By the  
definition of $b_\infty$, the point  $(b_\infty,0)$ which is in 
$D(a_\infty+\beta,b_\infty-\beta)$ is at a distance less than $2\eps$ 
from $\partial \Ome_n$ for big $n$ (see Figure \ref{dessin}). Then
there exists, for each big $n$, 
a point $q_n$ in $\partial\Ome_n$ such that $d_{\Ome_n} (q_n,
(b_\infty,0))\le 2\eps$. By \eqref{truc}, we can assume that for $n$
big enough: 
$$v_n(b_\infty,0)-v_n(O)\ge b_\infty-\eps$$
then :
\begin{align*}
v_n(O)&=v_n(O)-v_n(b_\infty,0)+v_n(b_\infty,0)\\
&\le \eps-b_\infty+ v_n(b_\infty,0)\\
&\le \eps-b_\infty+ 2\eps+\phi(q_n)=3\eps -b_\infty+ \phi(q_n)
\end{align*}
Besides
$$ v'_n(O)\ge \phi(q_n)-d_{\Ome_n}(O,q_n)\ge
\phi(q_n)-2\eps-b_\infty$$
Then $v_n(O)-v'_n(O)\le 5\eps$. The sequence $(p_n)$ is chosen in
$\Omet$ then $v(p_n)-v'(p_n)> a$ and then $v_n(O)-v'_n(O)\ge
a$. Then if $\eps$ is chosen such that $\eps<a/5$, we get a
contradiction and Lemma \ref{tech2} is proved.
\end{proof}

%% file: dessin.pstex_t
\begin{picture}(0,0)%
\includegraphics{dessin.pstex}%
\end{picture}%
\setlength{\unitlength}{4144sp}%
\begingroup\makeatletter\ifx\SetFigFont\undefined%
\gdef\SetFigFont#1#2#3#4#5{%
  \reset@font\fontsize{#1}{#2pt}%
  \fontfamily{#3}\fontseries{#4}\fontshape{#5}%
  \selectfont}%
\fi\endgroup%
\begin{picture}(7602,2352)(2239,-4573)
\put(5041,-3616){\makebox(0,0)[lb]{\smash{{\SetFigFont{14}{16.8}{\rmdefault}{\mddefault}{\updefault}{\color[rgb]{0,0,0}$0$}%
}}}}
\put(9136,-2401){\makebox(0,0)[lb]{\smash{{\SetFigFont{14}{16.8}{\rmdefault}{\mddefault}{\updefault}{\color[rgb]{0,0,0}$\partial\Ome_n$}%
}}}}
\put(9226,-4156){\makebox(0,0)[lb]{\smash{{\SetFigFont{14}{16.8}{\rmdefault}{\mddefault}{\updefault}{\color[rgb]{0,0,0}$q_n$}%
}}}}
\put(8551,-3481){\makebox(0,0)[lb]{\smash{{\SetFigFont{14}{16.8}{\rmdefault}{\mddefault}{\updefault}{\color[rgb]{0,0,0}$b_\infty$}%
}}}}
\put(2566,-3976){\makebox(0,0)[lb]{\smash{{\SetFigFont{14}{16.8}{\rmdefault}{\mddefault}{\updefault}{\color[rgb]{0,0,0}$D(a_\infty+\beta,b_\infty-\beta)$}%
}}}}
\end{picture}%